\documentclass[twoside,12pt]{article}
\usepackage{geometry}
\usepackage{amsmath,amssymb,amsfonts,amsthm}
\usepackage{txfonts}
\usepackage{todonotes}
\usepackage{tablists}
\usepackage{graphicx}
\usepackage{mathrsfs}
\usepackage{color}
\usepackage{float}
\usepackage{extarrows}
\usepackage{exscale}
\usepackage{titlesec}
\titleformat{\section}[block]{\large\center\sc}{\arabic{section}}{0.5em}{}[]
\usepackage{relsize}
\usepackage{color}
\definecolor{teal}{RGB}{67,225,128}
\usepackage{hyperref}
\hypersetup{hypertex=true,
            colorlinks=true,
            linkcolor=blue,
            anchorcolor=blue,
            citecolor=blue}
\usepackage[titles]{tocloft}
\usepackage{fancyhdr}
\theoremstyle{plain}
\newtheorem{theorem}{Theorem}[section]
\newtheorem{lemma}[theorem]{Lemma}

\newtheorem{proposition}[theorem]{Proposition}

\newtheorem{condition}[theorem]{Condition}

\newtheorem{remark}[theorem]{Remark}

\let\oldsection\section
\renewcommand\section{\setcounter{equation}{0}\oldsection}

\def\be{\begin{equation}}
\def\ee{\end{equation}}
\def\bes{\begin{equation*}}
\def\ees{\end{equation*}}
\def\bs{\begin{split}}
\def\es{\end{split}}
\def\bali{\begin{aligned}}
\def\eali{\end{aligned}}
\newcommand{\pf}{\noindent {\bf Proof. \hspace{2mm}}}
\def\bR{{\mathbb R}}
\def\un{\underbrace}
\def\al{\alpha}

\def\ve{\varepsilon}

\def\th{\theta}

\def\dl{\delta}
\def\Dl{\Delta}
\def\lt{\left}
\def\rt{\right}

\def\i{\infty}
\def\p{\partial}
\def\f{\frac}
\def\na{\nabla}
\def\o{\omega}
\def\O{\Omega}

\def\q{\quad}

\def\bl{\boldsymbol}

\def\mR{\mathbb{R}}
\def\mL{\mathcal{L}}

\allowdisplaybreaks

\def\cd{\cdot}
\def\les{\lesssim}

\begin{document}
\title{\bf\Large  A refined long time asymptotic bound for 3D axially symmetric Boussinesq system with zero thermal diffusivity}

\author{\normalsize\sc Zijin Li}

\date{}

\maketitle

\begin{abstract}
 In this paper, we obtain a refined temporal asymptotic upper bound of the global axially symmetric solution to the Boussinesq system with no thermal diffusivity. We show the spacial $W^{1,p}$-Sobolev ($2\leq p<\i$) norm of the velocity can only grow at most algebraically as $t\to+\i$. Under a signed potential condition imposed on the initial data, we further derive that the aforementioned norm is uniformly bounded at all times. Higher order estimates are also given: We find the $H^1$ norm of the temperature fluctuation grows sub-exponentially as $t\to+\i$. Meanwhile, for any $m\geq1$, we deduce that the $H^m$-temporal growth of the solution is slower than a double exponential function. As a result, these improve the results in \cite{HR:2010AIHP} where the authors only provided rough temporal asymptotic upper bounds while proving the global well-posedness.

\medskip



\end{abstract}

\tableofcontents

\section{Introduction}\label{SEC1}
The Boussinesq system is commonly used in describing the motion of ocean or atmospheric dynamics. It is derived from the density dependent incompressible Navier-Stokes equations by using the Boussinesq approximation by neglecting the density dependence in all the terms besides the one involving gravity. It reads that
\be\label{Bous}
\left\{\begin{array}{l}
\partial_t v+v \cdot \nabla v-\mu\Delta v+\nabla P=\rho \bl{e_3}, \quad(t, x) \in \mathbb{R}_{+} \times \mathbb{R}^3, \\
\partial_t \rho+v \cdot \nabla \rho=0, \\
\operatorname{div} v=0, \\
v(0,x)=v_0, \quad \rho(0,x)=\rho_0 .
\end{array}\right.
\ee
Here, the velocity $v=\left(v_1, v_2, v_3\right)$ is a divergence-free three dimensional vector field, while $P\in\bR$ and $\rho\in\bR$ represent the pressure and the temperature fluctuation, respectively. $\bl{e_3}=(0,0,1)^T$ is the unit vector in the vertical direction. The coefficient $\mu>0$ is the Reynolds number that measures the strength of heat conductivity, and in the following we assume it equals to one without loss of generality. In terms of physics, equation \eqref{Bous}$_1$ illustrates the conservation law of the momentum under the influence of the buoyant effect $\rho \bl{e_3}$. Equation \eqref{Bous}$_2$ describes temperature fluctuation with diffusivity, while equation \eqref{Bous}$_3$ describes fluid incompressibility.

Notice that when the initial density $\rho_0$ is identically zero or constant, then the system \eqref{Bous} reduces to the classical incompressible Navier-Stokes equation:
\[
\left\{\begin{array}{l}
\partial_t v+v \cdot \nabla v-\Delta v+\nabla P=0, \\
\operatorname{div} v=0, \\
v(0,x)=v_0 .
\end{array}\right.
\]
Therefore, one cannot expect to have a better theory for the Boussinesq system than for the Navier-Stokes equations. For the two-dimensional case, the global well-posedness of Boussinesq systems has received a great deal of attention in the past few decades. Global well-posedness has been shown in various function spaces and for different type of viscosities, we refer for example to \cite{AH:2007JDE, Brenier2009, CaoWu2013, Chae:2006ADV, Danchin2008, HK:2007ADE, HK:2009INDI, HL:2005DCDS, KW2020JDDE, Larios2013}. For full three-dimensional systems, even when $\rho$ is missing, the global existence of the strong solution with general initial data is far from being solved. Danchin-Paicu \cite{Danchin2009} showed a global well-posedness result for small initial data belonging to some critical Lorentz spaces. Abidi et al. \cite{AHK:2011DCDS} and Hmidi-Rousset \cite{HR:2010AIHP,HR:2011JFA} proved the global well-posedness of the Cauchy problem for the 3D axisymmetric Boussinesq system without swirl. When the swirl component presents, recently Li-Pan \cite{Li-Pan2022DCDS} derived a single-component Prodi-Serrin-type regularity criterion for strong solutions. They also gave a single-component Beale-Kato-Majda-type regularity criterion for related inviscid equations which has positive thermal diffusivity in \cite{Li-Pan2023ZAMP}.

\subsection{The axisymmetric Boussinesq system in the cylindrical coordinates}
Our proof will be partially presented out in the cylindrical coordinates $(r, \theta, z)$. That is, for $x=\left(x_1, x_2, x_3\right) \in \mathbb{R}^3$:
$$
r=\sqrt{x_1^2+x_2^2}, \quad \theta=\arctan \frac{x_2}{x_1}, \quad z=x_3 .
$$
And the axisymmetric solution of system \eqref{Bous} is given by
$$
\left\{
\begin{aligned}
v&=v_r(t, r, z) \bl{e_r}+v_\theta(t, r, z) \bl{e_\theta}+v_z(t, r, z) \bl{e_z} \\
\rho&=\rho(t, r, z)
\end{aligned}
\right.
$$
where the basis vectors $\bl{e_r}, \bl{e_\theta}, \bl{e_z}$ are
$$
\bl{e_r}=\left(\frac{x_1}{r}, \frac{x_2}{r}, 0\right), \quad \bl{e_\theta}=\left(-\frac{x_2}{r}, \frac{x_1}{r}, 0\right), \quad \bl{e_z}=(0,0,1)
$$
After a direct but tedious calculation, one derives the axially symmetric Boussinesq system in the cylindrical coordinates as follows:
\be\label{(2.1)}
\left\{
\begin{split}
&\partial_t v_r+\left(v_r \partial_r+v_z \partial_z\right) v_r-\frac{\left(v_\theta\right)^2}{r}+\partial_r P=\left(\Delta-\frac{1}{r^2}\right) v_r, \\[1mm]
&\partial_t v_\theta+\left(v_r \partial_r+v_z \partial_z\right) v_\theta+\frac{v_\theta v_r}{r}=\left(\Delta-\frac{1}{r^2}\right) v_\theta, \\[1mm]
&\partial_t v_z+\left(v_r \partial_r+v_z \partial_z\right) v_z+\partial_z P=\Delta v_z+\rho, \\[1mm]
&\p_t\rho+\left(v_r \partial_r+v_z \partial_z\right)\rho=0,\\[1mm]
&\nabla \cdot v=\partial_r v_r+\frac{v_r}{r}+\partial_z v_z=0,
\end{split}
\right.
\ee
where
\[
\Delta=\frac{\partial^2}{\partial r^2}+\frac{1}{r} \frac{\partial}{\partial r}+\frac{\partial^2}{\partial z^2}
\]
is the usual Laplacian operator. As a result of the uniqueness of local solutions, it is clear that if the initial data satisfy $v_0\cdot\bl{e_\th}=0$, then the solution of \eqref{(2.1)} will be the following:
\[
v=v_r(t,r,z)\bl{e_r}+v_z(t,r,z)\bl{e_z},\q\rho=\rho(t,r,z).
\]
And \eqref{(2.1)} can be simplified as
\be\label{2.2}
\left\{
\begin{split}
&\partial_t v_r+\left(v_r \partial_r+v_z \partial_z\right) v_r+\partial_r P=\left(\Delta-\frac{1}{r^2}\right) v_r, \\[1mm]
&\partial_t v_z+\left(v_r \partial_r+v_z \partial_z\right) v_z+\partial_z P=\Delta v_z+\rho, \\[1mm]
&\p_t\rho+\left(v_r \partial_r+v_z \partial_z\right)\rho=0,\\[1mm]
&\nabla \cdot v=\partial_r v_r+\frac{v_r}{r}+\partial_z v_z=0.
\end{split}
\right.
\ee
Taking the curl of the first and the second equations in \eqref{2.2}, one derives $\o_\th:=\p_zv_r-\p_rv_z$ satisfies
\be\label{omega}
\partial_t \omega_\th+\left(v_r\p_r+v_z\p_z\right)\omega_\th=\left(\Dl-\f{1}{r^2}\right)\o_\th+\frac{v_r}{r} \omega_\th-\p_r\rho .
\ee
This implies $\Omega:=\frac{\o_\th}{r}$ enjoys that
\[
\partial_t \Omega+v \cdot \nabla \Omega=\left(\Delta+\frac{2}{r} \partial_r\right) \Omega-\frac{\partial_r \rho}{r}.
\]

We finish this subsection with a list of notations that will appear throughout the paper.
\begin{itemize}
\item $C_{a,b,...}$ denotes a positive constant depending on $a,\,b,\,...$ which may be different from line to line. Likewise, we use $C_{0,...}$ to denote a constant that also depends on initial data.

\item $A\lesssim B$ means $A\leq CB$. Meanwhile, $A\simeq B$ means both $A\lesssim B$ and $B\lesssim A$.

\item $[\mathcal{A},\,\mathcal{B}]=\mathcal{A}\mathcal{B}-\mathcal{B}\mathcal{A}$ denotes the commutator of the operator $\mathcal{A}$ and the operator $\mathcal{B}$.

\item $\mathfrak{H}$ stands for a multi-index such that $\mathfrak{H}=(h_1,h_2,h_3)$ where $h_1,h_2,h_3\in\mathbb{N}\cup\{0\}$ and $|\mathfrak{H}|=h_1+h_2+h_3$, $\nabla^\mathfrak{H}=\p_{x_1}^{h_1}\p_{x_2}^{h_2}\p_{x_3}^{h_3}$.

\item We use standard notations for Lebesgue and Sobolev functional spaces in $\mathbb{R}^3$: For $1\leq p\leq\infty$ and $k\in\mathbb{N}$, $L^p$ denotes the Lebesgue space with norm
\[
\|f\|_{L^p}:=
\lt\{
\begin{aligned}
&\left(\int_{\mathbb{R}^3}|f(x)|^pdx\right)^{1/p},\quad 1\leq p<\infty,\\
&\mathop{ess sup}_{x\in\mathbb{R}^3}|f(x)|,\quad\quad\quad\quad p=\infty.\\
\end{aligned}
\rt.
\]

\item $W^{k,p}$ denotes the usual Sobolev space with its norm
\[
\begin{split}
\|f\|_{W^{k,p}}:=&\sum_{0\leq|\mathfrak{H}|\leq k}\|\nabla^\mathfrak{H} f\|_{L^p}\,.\\
\end{split}
\]
We also simply denote $H^k$ and $\dot{H}^k$ instead of $W^{k,p}$ and $\dot{W}^{k,p}$ provided $p=2$.
\item For a function $f\in L^p\cap L^q$ with $1\leq p,q\leq \infty$, we denote its Yudovich-type norm as
\[
\|f\|_{L^p\cap L^q}=\max\left\{\|f\|_{L^p},\|f\|_{L^q}\right\}.
\]
\item For any Banach space $X$, we say $v:\,[0,T]\times\mathbb{R}^3\to\mathbb{R}$ belongs to the Bochner space $L^q_TX$, if
\[
\|v(t,\cdot)\|_{X}\in L^q\big((0,T)\big).
\]
\item For a fixed real number $\varsigma$, we denote $\varsigma_+$ a number larger than but close to $\varsigma$.

\end{itemize}

\subsection{Main results}
Now we are ready for main results of the paper. To state it, we first introduce the following condition on the potential energy:
\begin{condition}[bounded \& signed potential]\label{COND1}
We say the axisymmetric initial data $(v_0,\rho_0)$ satisfies the bounded \& signed potential condition if
\begin{itemize}
\item[(i).]
$(v_r)_0$ is even, and $(v_z)_0$ and $\rho_0$ are odd symmetric in $ x_3 $;
\item[(ii).]
$\rho_0\leq 0$, for all $x_3>0$;
\item[(iii).]
It holds that
\[
\int_{\mR^3}\rho_0x_3dx>-\i.
\]
\end{itemize}
\end{condition}

\qed

By uniqueness, the symmetry property \emph{(i)} above will keep for any $t>0$. This indicates $v_z\equiv 0$ on $z=0$. With the help of the trajectory mapping argument, one has
\be\label{EG0}
\rho(t,x)\leq 0,\q\text{for all}\q (t,x)\in\mathbb{R}_+\times\mathbb{R}^2\times\mathbb{R}_+\,.
\ee
Here goes the reason: Denote $X(t,\cdot):\,\mathbb{R}^3\to\mathbb{R}^3$ the particle trajectory mapping of the velocity $v$, which solves the initial value problem:
\bes
\f{\p X(t,\zeta)}{\p t}=v(t,X(t,\zeta)),\quad X(0,\zeta)=\zeta.
\ees
Integrating \eqref{2.2}$_3$ along the particle trajectory mapping, we have
\[
\rho(t,X(t,\zeta))=\rho_0(\zeta).
\]
Since $v_z(t,r,0)\equiv0$, one observes a trajectory $X(t,\zeta)$ could not go across the hyperplane $\{z=0\}$, which concludes the validity of \eqref{EG0}. Therefore, under the Condition \ref{COND1} \emph{(i)\,--\,(ii)}, one has the potential energy
\be\label{PT}
Po(t):=\int_{\mathbb{R}^3}\rho(t,x)x_3dx\leq0,\q\forall t\geq 0.
\ee
Item \emph{(iii)} of Condition \ref{COND1} is rather physical, which indicates the finiteness of the initial potential energy $Po(t)$.

\begin{remark}
Condition \ref{COND1} in this paper is motivated by the assumption (A2) in \cite{KPY2022} where the authors introduced there to do research on lower temporal bound of 2D Boussinesq equations and 3D axially symmetric Euler equations. Following their idea, we can derive the temperature fluctuation of the system \eqref{2.2} enjoys the following temporal lower bound:
\[
\|\rho(t,\cd)\|_{\dot{H}^s}\geq Ct^{s/8}.
\]
Consider the case of $s=1$, this lower bound is algebraical as $t\to\i$, while the upper bound which will be obtained in Theorem \ref{Main2} is sub-exponential. There is still a clear gap between them. Currently, in the author's opinion, the strategy of this work does not appear to be strong enough to fill it.
\end{remark}

\qed

Thanks to the global well-posedness for the Boussinesq system with axisymmetric data obtained in \cite{HR:2010AIHP}, in this paper we only focus on the refined long time asymptotic behavior for global smooth solutions to system \eqref{2.2}. Our first result shows the ${W}^{1,p}$ norm of $v(t,\cd)$ is uniformly bounded with $t>0$ if Condition \ref{COND1} is satisfied.

\begin{theorem}\label{Main1}
Let $(v,\rho)$ be a smooth global axially symmetric solution of \eqref{Bous}, with its divergence-free initial data $v_0\in H^1$, $v_0\cdot\bl{e_\th}=0$, $\o_0\in L^\i$, together with $\f{(\o_\th)_0}{r}\in L^2$ and $\rho_0\in L^2\cap L^\i$. If the initial data $(v_0,\rho_0)$ further satisfy Condition \ref{COND1}, then there exists $C_{0,p}>0$ that
\[
\sup_{0\leq t<\i}\|v(t,\cdot)\|_{W^{1,p}}\leq C_{0,p}.
\]
Here $2\leq p<\i$.
\end{theorem}

\qed

In Theorem \ref{Main1}, if Condition \ref{COND1} no longer holds for the initial data, one can still have the $L^\i_t{W}^{1,p}\cap L^q_t\dot{W}^{2,p}$ norm of $v$ grows at most algebraically as $t\to\i$. This improves the known result in \cite{HR:2010AIHP} where the authors proved at most exponential growth. Here is the result:

\begin{proposition}\label{COR1}
Let $(v,\rho)$ be a smooth global axially symmetric solution of \eqref{Bous}, with its divergence-free initial data $v_0\in H^1$, $v_0\cdot\bl{e_\th}=0$, $\o_0\in L^\i$, together with $\f{(\o_\th)_0}{r}\in L^2$ and $\rho_0\in L^2\cap L^\i$. Then there exists $M>0$, depending only on $p$ and $q$, such that:
\[
\sup_{0\leq s\leq t}\|v(s,\cdot)\|_{W^{1,p}}+\left(\int_{0}^t\|\na^2 v(s,\cd)\|_{L^p}^qds\right)^{1/q}\leq C_{0,p,q}(1+t)^M.
\]
Here $2\leq p<\i$, $1\leq q<\i$\,.
\end{proposition}

\qed

The next result gives higher-order energy estimates of both $v$ and $\rho$ under the \emph{bounded \& signed potential} condition:

\begin{theorem}\label{Main2}
Let $(v,\rho)$ be a smooth global axially symmetric solution of \eqref{Bous}, with its initial data $(v_0,\rho_0)\in H^m\times H^m$ ($3\leq m\in\mathbb{N}$) , and $\na\cd v_0=v_0\cdot\bl{e_\th}=0$. If $(v_0,\rho_0)$ further satisfies Condition \ref{COND1}, then for any $\dl>0$, we have
\be\label{EH1}
\begin{split}
\left\|\na \rho(t,\cdot)\right\|_{L^p}\les_{0,p,\dl}&\exp\left(t^{\f{4}{5}+\dl}\right),\quad \forall t\in[0,\i),\q 1\leq p\leq\i\,,
\end{split}
\ee
and
\be\label{EH2}
\begin{split}
\left\|\na^m (v,\rho)(t,\cdot)\right\|_{L^2}^2\les_{0,m,\dl}&\exp\left(\exp\left(t^{\f{4}{5}+\dl}\right)\right),\quad \forall t\in[0,\i)\,.
\end{split}
\ee
\end{theorem}

\qed

In Theorem \ref{Main2}, if the initial data do not satisfy Condition \ref{COND1}, one can still have rather weaker estimates by changing the exponent $(4/5)_+$ in \eqref{EH1} and \eqref{EH2} to a certain positive constant $M$. Here is the result:

\begin{proposition}\label{Main3}
Let $(v,\rho)$ be a smooth global axially symmetric solution of \eqref{Bous}, with its initial data $(v_0,\rho_0)\in H^m\times H^m$ ($3\leq m\in\mathbb{N}$) , and $\na\cd v_0=v_0\cdot\bl{e_\th}=0$. Then
\[
\begin{split}
\left\|\na \rho(t,\cdot)\right\|_{L^p}&\les_{0,p}\exp\left(t^{\tilde{M}}\right),\q 1\leq p\leq\i\,;\\[1mm]
\left\|\na^m (v,\rho)(t,\cdot)\right\|_{L^2}^2&\les_{0,m}\exp\left(\exp\left(t^{\tilde{M}}\right)\right)\,,
\end{split}
\]
for all $t\in[0,\i)$. Here $\tilde{M}$ is a positive constant that is independent with $p$ or $m$.
\end{proposition}

\qed

\subsection{Strategy of the proof}
Now we outline the idea in proving main results of the paper. With the help of the \emph{bounded \& signed potential condition}, one derives the following fundamental energy estimate:
\be\label{OOO}
\|\rho(t,\cd)\|_{L^p}+\|v(t,\cd)\|_{L^2}^2+\int_0^t\|\nabla v(s,\cd)\|_{L^2}^2 d s \leq C_0
\ee
for any $p\in[2,\i]$. Unlike deriving the long time behavior of solutions to the 2D Boussinesq equations in \cite{KW2020JDDE}, we must overcome difficulties caused by the vortex stretching effect in the 3D problem. Motivated by Hmidi-Rousset \cite{HR:2010AIHP}, we prove a self-closed energy estimate of the quantity $\O-\left(\Delta+\frac{2}{r} \partial_{r}\right)^{-1} \frac{\partial_{r}\rho}{r}$, and then it follows that:
\be\label{OO}
\|\O(t,\cd)\|_{L^2}\leq C_0\,.
\ee
When trying to obtain this uniform-in-time bound of $\O$, one must be very careful to avoid applying the Gr\"onwall inequality with an exponential factor. To this end, we introduce the following Hardy-type inequality
\[
\left\|\f{f}{\sqrt{r}}\right\|_{L^3}\lesssim\|\na f\|_{L^2},\q\text{for}\q f\in C_c^\i(\mathbb{R}^3)
\]
and take a completely different path when dealing with the major term. With the help of \eqref{OOO} and \eqref{OO}, we are able to show
\[
\left\|\o_\theta(t,\cd)\right\|_{L^{2^n}}^{2^n}+\int_0^t\big\|\nabla \o_\theta^{2^{n-1}}(s,\cd)\big\|_{L^2}^2ds+\int_0^t\Big\|\frac{\o_\theta^{2^{n-1}}}{r}(s,\cd)\Big\|_{L^2}^2ds\leq C_{0,n}\,,\q\text{for}\q n\in\mathbb{N}\cap\{n\geq 2\}
\]
by induction. This concludes Theorem \ref{Main1}. Using the maximal regularity of the heat flow, we are able to deduce
\be\label{O}
\int_0^t\|\na^2v(s,\cd)\|^p_{L^{2^n}}ds\leq C_{0,p,n}t\,,
\ee
and then we arrive at
\[
\int_0^{t}\|\nabla v(s,\cdot)\|_{L^\infty}ds\lesssim t^{\left(4/5\right)_+}
\]
by an interpolation between \eqref{OOO} and \eqref{O}. Thus the higher-order estimates in Theorem \ref{Main2} are derived by a routine energy estimate and an application of the Gr\"onwall inequality.

If Condition \ref{COND1} is ignored, instead of \eqref{OOO}, one can only derive
\[
\|v(t,\cd)\|_{L^2}^2+\int_0^t\|\nabla v(s,\cd)\|_{L^2}^2 d s \leq C_0(1+t)^2\,.
\]
This results in the algebraical temporal growth of quantities of $\O$ and $\o_\th$, which lead to results in Proposition \ref{COR1} and Proposition \ref{Main3}. With a careful calculation, one can find the optimal growth order therein, but we will not pursue it in this paper.

The rest of this paper is organized as follows. In Section \ref{PRE}, we provide some useful Lemmas concerning interpolation inequalities, some $L^p$ boundedness of singular operators related to the problem, a commutator estimate by Kato-Ponce, a Hardy type inequality, and the maximal regularity for the heat flow. Finally, main results will be proved in Section \ref{Main}.

\section{Preliminary}\label{PRE}
Some well-known lemmas will be listed in this section without detailed proof. At the beginning, let us introduce the well-known $Gagliardo-Nirenberg$ interpolation inequality.
\begin{lemma}[Gagliardo-Nirenberg]\label{LEMGN}
Given $q,r\in[1,\i]$ and $j,m\in\mathbb{N}\cup\{0\}$ with $j\leq m$. Suppose that $f\in L^q(\mathbb{R}^d)\cap\dot{W}^{m,r}(\mathbb{R}^d)$ and there exists a real number $\al\in[j/m,1]$ such that
\[
\frac{1}{p}=\frac{j}{d}+\al\left(\frac{1}{r}-\frac{m}{d}\right)+\frac{1-\al}{q}.
\]
Then $f\in\dot{W}^{j,p}(\mathbb{R}^d)$ and there exists a constant $C>0$ such that
\[
\|\na^jf\|_{L^p(\mathbb{R}^d)}\leq C\|\na^m f\|^\al_{L^r(\mathbb{R}^d)}\|f\|^{1-\al}_{L^q(\mathbb{R}^d)},
\]
except the following two cases:
\begin{itemize}
\item[I.] $j=0$, $mr<d$ and $q=\infty$; (In this case it is necessary to assume also that either $|u|\to 0$ at infinity, or $u\in L^s(\mathbb{R}^d)$ for some $s<\infty$.)

\item[II.] $1<r<\infty$ and $m-j-d/r\in\mathbb{N}$. (In this case it is necessary to assume also that $\alpha<1$.)
\end{itemize}
\end{lemma}

\qed

The following lemma was introduced by Hmidi-Rousset \cite{HR:2010AIHP} where the authors derived regularity of the axisymmetric Boussinesq system without swirl. It states the $L^p$-boundedness of two operators related to axially symmetric vector fields.

\begin{lemma}\label{LEMET11}
Denote $\mathcal{L}=\left(\Delta+\frac{2}{r} \partial_{r}\right)^{-1} \frac{\partial_{r}}{r}$ and $\tilde{\mathcal{L}}=\left(\Delta+\frac{2}{r} \partial_{r}\right)^{-1} \frac{\partial_{z}}{r} .$ Suppose $\rho \in H^{2}\left(\mathbb{R}^{3}\right)$ be axisymmetric, then for every $p \in[2,+\infty),$ there
exists an absolute constant $C_{p}>0$ such that
\[
\|\mathcal{L} \rho\|_{L^{p}} \leq C_{p}\|\rho\|_{L^{p}}, \quad\|\tilde{\mathcal{L}} \rho\|_{L^{p}} \leq C_{p}\|\rho\|_{L^{p}}.
\]
Moreover, for any smooth axisymmetric function $f,$ we have the identity
\be\label{LEMET11L}
\mathcal{L} \partial_{r} f=\frac{f}{r}-\mathcal{L}\left(\frac{f}{r}\right)-\partial_{z} \tilde{\mathcal{L}} f.
\ee
\end{lemma}

\pf The detailed proof can be found in Proposition 3.1, 3.2 and Lemma 3.3 in \cite{HR:2010AIHP}. We omit the details here.

\qed

The following famous commutator estimate will be applied later in our proof.

\begin{lemma}\label{LEMET1}
Let $m\in\mathbb{N}$, $m\geq 2$, and $f,g,k\in C^\infty_0(\mathbb{R}^3)$. Then the following estimate holds:
\be\label{E1}
\begin{split}
\left|\int_{\mathbb{R}^3}[\nabla^m,\,f\cdot\nabla]g\nabla^m kdx\right|\leq&\,C\,\left\|\nabla^{m}(f,g,k)\right\|_{L^2}^2\|\nabla (f,\,g)\|_{L^\infty}.
\end{split}
\ee
\end{lemma}
\pf Applying H\"{o}lder's inequality, one derives
\be\label{2.77777}
\left|\int_{\mathbb{R}^3}[\nabla^m,\,f\cdot\nabla]g\nabla^m kdx\right|\leq\|[\nabla^m,\,f\cdot\nabla]g\|_{L^{2}}\|\nabla^mk\|_{L^2}.
\ee
Due to the commutator estimate by Kato-Ponce \cite{Kato1988}, it follows that
\be\label{2.88888}
\|[\nabla^m,\,f\cdot\nabla]g\|_{L^{2}}\leq C\left(\|\nabla f\|_{L^\infty}\|\nabla^mg\|_{L^2}+\|\nabla g\|_{L^\infty}\|\nabla^m f\|_{L^2}\right).
\ee
Then \eqref{E1} follows from plugging \eqref{2.88888} into \eqref{2.77777}.

\qed

Next we give a Sobolev-Hardy inequality. We omit the detailed proof since it could be found in the Lemma 2.4 of \cite{ChenFangZhang2017}.
\begin{lemma}\label{CFZ}
Set $\mathbb{R}^{d}=\mathbb{R}^{k} \times \mathbb{R}^{d-k}$ with $2 \leq k \leq d,$ and write $x=\left(x^{\prime}, z\right) \in \mathbb{R}^{k} \times \mathbb{R}^{d-k}$. For
$1<q<d, 0 \leq \theta \leq q$ and $\theta<k,$ let $q_{*} \in\left[q, \frac{q(d-\theta)}{d-q}\right]$. Then there exists a positive constant $C=C(\theta, q, d, k)$ such that for all $f \in C_{0}^{\infty}\left(\mathbb{R}^{d}\right),$
\[
\left(\int_{\mathbb{R}^{d}} \frac{|f|^{q_{*}}}{\left|x^{\prime}\right|^{\theta}} d x\right)^{\f{1}{q_*}} \leq C\|f\|_{L^q}^{\frac{d-\theta}{q_{*}}-\frac{d}{q}+1}\|\nabla f\|_{L^q}^{\frac{d}{q}-\frac{d-\theta}{q_*}}.
\]
In particular, we pick $d=3$, $k=2$, $q_*=3$, $q=2$, $\th=3/2$ and assume $r=$ $\sqrt{x_{1}^{2}+x_{2}^{2}} .$ Then there exists a positive constant $C$ such that for all $f \in C_{0}^{\infty}\left(\mathbb{R}^{d}\right)$
\be\label{HEE}
\left\|\f{f}{\sqrt{r}}\right\|_{L^3}\leq C\|\na f\|_{L^2}\,.
\ee
\end{lemma}
\qed



Now we recall the standard maximal regularity of heat flows in $L^r_TL^p$-type spaces. Detailed proof could be found in \cite[Theorem 7.3]{Rieusset2002} for instance.
\begin{lemma}[Maximal $L^r_TL^p$ regularity for the heat flow]\label{MRP}
Let us define the operator $\mathcal{A}$ by the formula
\[
\mathcal{A}: \quad f \longmapsto \int_{0}^{t} \nabla^{2} e^{(t-s) \Delta} f(s, \cdot) d s.
\]
Then $\mathcal{A} \text { is bounded from } L^{r}\left(0, T; L^{p}(\mathbb{R}^{d})\right) \text { to } L^{r}\left( 0, T; L^{p}(\mathbb{R}^{d})\right)$ for every  $T\in(0,\infty]$ and
$1<p, r<\infty$. Moreover, there holds:
\[
\left\|\mathcal{A} f\right\|_{ L_{T}^{r}L^{p}} \leq C\|f\|_{ L_{T}^{r}L^{p}}.
\]
\end{lemma}
\qed

\section{Proof of main results}\label{Main}

\subsection{Fundamental energy bound}\label{SEU}
The first step states the fundamental energy estimate of the solution:
\begin{proposition}
Let $(v,\rho)$ be a global smooth solution of \eqref{Bous}, then we have:
\begin{itemize}
\item[(1).] Let $\rho_0\in L^p$ for $p \in[1, \infty]$, and $t \in \mathbb{R}_{+}$, we have
\be\label{ERHO}
\|\rho(t,\cd)\|_{L^p} \leq\left\|\rho_0\right\|_{L^p} ;
\ee
\item[(2).] For $(v_0,\rho_0)$ satisfies Condition \ref{COND1} such that $v_0\in L^2$, then for all $t \in \mathbb{R}_{+}$, we have
\be\label{Ev}
\|v(t,\cd)\|_{L^2}^2+\int_0^t\|\nabla v(s,\cd)\|^2 d s \leq C_0\,.
\ee
\end{itemize}
Here $C_0$ is a universal constant depends only on the initial data $(v_0,\rho_0)$\,.
\end{proposition}
\pf For $1\leq p<\i$, estimate \eqref{ERHO} is classical for transport equations with a divergence-free vector field. And when $p=\i$, \eqref{ERHO} is known as the maximum principle. Now we derive \eqref{Ev} by carrying out the standard $L^2$ energy estimate of the system \eqref{Bous}. Multiplying $v$ on both sides of \eqref{Bous}$_1$ and integrating over $\mathbb{R}^3$, we find
\be\label{FE1}
\f{1}{2}\f{d}{dt}\|v(t,\cdot)\|_{L^2}^2+\|\na v(t,\cdot)\|_{L^2}^2=\int_{\mathbb{R}^3}\rho v_3dx.
\ee
Meanwhile, multiplying $x_3$ on both sides of \eqref{Bous}$_2$ and integrating over $\mR^3$, one arrives
\be\label{FE2}
\f{d}{dt}\int_{\mR^3}\rho x_3dx=-\int_{\mR^3}x_3\na\cd(\rho v)dx=\int_{\mR^3}\rho v_3dx.
\ee
Subtracting \eqref{FE2} from \eqref{FE1}, and integrating with the temporal variable over $(0,t)$, one deduces
\[
\f{1}{2}\|v(t,\cdot)\|_{L^2}^2-\int_{\mR^3}\rho(t,x)x_3dx+\int_0^t\|\nabla v(s,\cdot)\|_{L^2}^2ds=\f{1}{2}\|v_0\|_{L^2}^2-\int_{\mR^3}\rho_0(x)x_3dx.
\]
Using \eqref{PT}, one concludes the fundamental energy estimate of $\eqref{Bous}$ that:
\be\label{Fund}
\f{1}{2}\|v(t,\cdot)\|_{L^2}^2+\int_0^t\|\nabla v(s,\cdot)\|_{L^2}^2ds\leq C_0,\q\forall t\geq 0\,.
\ee
This completes the proof of the proposition.

\qed

\subsection{Estimate of a reformulated equation }
Unlike the 2D Boussinesq system, the presence of the vortex stretching term in the 3D vorticity equation poses significant challenges:
\[
\partial_t \omega_\th+\left(v_r\p_r+v_z\p_z\right)\omega_\th=\left(\Dl-\f{1}{r^2}\right)\o_\th+\un{\,\,\,\,\frac{v_r}{r} \omega_\th\,\,\,\,}_{\text{vortex stretching term}}-\p_r\rho .
\]
To overcome this difficulty, we first derive an $L^2$ a priori bound of $\O=\f{\o_\th}{r}$. To achieve a refined long time asymptotic behavior of the solution, one must take great care to avoid the temporal growth of $\|\O(t,\cd)\|_{L^2}$ at this point in the analysis.
\begin{proposition}\label{PROP3.2}
Let $(v,\rho)$ be a global smooth solution of axially symmetric Boussinesq system \eqref{2.2}, with its initial data $(v_0,\rho_0)$ satisfying Condition \ref{COND1}, and $\rho_0\in L^2\cap L^3$, $v_0\in L^2$, $\na\cdot v_0=v_0\cdot\bl{e_\th}=0$, $\f{(\o_\th)_0}{r}\in L^2$. Then $\O=\f{\o_\th}{r}$ satisfies
\[
\|\O(t,\cd)\|_{L^2}\leq C_0\,.
\]
\end{proposition}
\pf Recalling Lemma \ref{LEMET11}, acting $\mathcal{L}$ on \eqref{2.2}$_3$, we get
\be\label{211}
\partial_t \mathcal{L} \rho+v \cdot \nabla \mathcal{L} \rho=-[\mathcal{L}, v \cdot \nabla] \rho\,.
\ee
Denoting $L:=\Omega-\mathcal{L} \rho$, by the equation of $\O$:
\be\label{Omega1}
\partial_t \Omega+v \cdot \nabla \Omega=\left(\Delta+\frac{2}{r} \partial_r\right) \Omega-\frac{\partial_r \rho}{r}\,,
\ee
and noticing that $\f{\p_r\rho}{r}=\left(\Delta+\frac{2}{r} \partial_r\right)\mathcal{L}\rho\,$, we have the following reformulated equation by subtracting \eqref{211} from \eqref{Omega1}:
\be\label{EL}
\partial_t L+v \cdot \nabla L=\left(\Delta+\frac{2}{r} \partial_r\right) L+[\mathcal{L}, v \cdot \nabla] \rho\,.
\ee
Taking the $L^2$ inner product of \eqref{EL} with $L$ and using integration by parts, noticing that
\[
\int_{\mathbb{R}^3}\f{2}{r}\p_rL\cdot L dx=2\pi\int_{\mathbb{R}}\int_0^\infty\f{1}{r}\p_rL^2rdrdz=-2\pi\int_{\mathbb{R}}|L(t,0,z)|^2dz\leq 0\,,
\]
one arrives at
\be\label{E01}
\begin{aligned}
\frac{1}{2} \frac{d}{d t}\|L(t,\cd)\|_{L^2}^2+\|\nabla L(t,\cd)\|_{L^2}^2\leq & \int_{\mathbb{R}^3} \mathcal{L}(v \cdot \nabla \rho) L d x-\int_{\mathbb{R}^3} v \cdot \nabla(\mathcal{L} \rho) L d x\\
= & \un{\int_{\mathbb{R}^3} \mathcal{L}(v \cdot \nabla \rho) L d x}_{I_1}+\un{\int_{\mathbb{R}^3}(\mathcal{L} \rho) v \cdot \nabla L d x}_{I_2}\,.\\
\end{aligned}
\ee
Here in $I_1$,
\[
\begin{aligned}
\mathcal{L}(v \cdot \nabla \rho) &=\mathcal{L}(\nabla \cdot(\rho v)) \\
&=\mathcal{L}\left(\partial_r\left(\rho v_r \right)\right)+\mathcal{L}\left(\frac{\rho v_r}{r} \right)+\mathcal{L}\left(\partial_z\left(\rho v_z \right)\right)\\
&=\frac{\rho v_r}{r}-\partial_z \tilde{\mathcal{L}}\left(\rho v_r\right)+\partial_z \mathcal{L}\left( \rho v_z\right).
\end{aligned}
\]
Here the third equality follows from \eqref{LEMET11L} in Lemma \ref{LEMET11} and the fact that $\p_z$ commutes with $\mathcal{L}$. Using integration by parts, we arrive
\[
I_1=\un{\int_{\mathbb{R}^3} \frac{v_r}{r} \rho L d x}_{J_1}+\un{\int_{\mathbb{R}^3} \tilde{\mathcal{L}}\left(\rho v_r\right) \partial_z L d x}_{J_2}-\un{\int_{\mathbb{R}^3} \mathcal{L}\left(\rho v_z\right) \partial_z L d x}_{J_3}.
\]
For the major term $J_1$, we apply H\"older's inequality to derive
\[
\begin{split}
|J_1|&\leq\Big\|\f{v_r}{\sqrt{r}}(t,\cd)\Big\|_{L^3}\|\rho(t,\cd)\|_{L^3}\Big\|\f{L}{\sqrt{r}}(t,\cd)\Big\|_{L^3}\,.
\end{split}
\]
Using the Hardy-type inequality \eqref{HEE}, one arrives at
\[
|J_1|\les\|\na v(t,\cd)\|_{L^2}\|\rho(t,\cd)\|_{L^3}\|\na L(t,\cd)\|_{L^2}\,.
\]
Then, by Young's inequality and the fundamental energy estimate \eqref{ERHO}, one concludes that
\[
|J_1|\leq\f{1}{4}\|\na L(t,\cd)\|_{L^2}^2+C\|\rho_0\|_{L^3}^2\|\na v(t,\cd)\|_{L^2}^2\,.
\]
\begin{remark}
Usually, one estimates $J_1$ in the following way:
\[
\begin{aligned}
\left|J_1\right| & \les\Big\|\na\frac{v_r}{r}(t,\cd)\Big\|_{L^2}\|\rho(t,\cd)\|_{L^3}\|L(t,\cd)\|_{L^2} \\
& \leq\left(\|L(t,\cd)\|_{L^2}+\|\mathcal{L} \rho(t,\cd)\|_{L^2}\right)\|\rho(t,\cd)\|_{L^3}\|L(t,\cd)\|_{L^2}\,.
\end{aligned}
\]
Here one has applied the following famous estimate (see for detailed proof in \cite{Lz:2015JDE} (equation (A.5)) and  \cite{Miao-Zheng2013} (Proposition 2.5))
\[
\Big\|\na\f{v_r}{{r}}(t,\cd)\Big\|_{L^p}\leq C_p\|\O(t,\cd)\|_{L^p},\q\text{for } 1<p<+\i\,,
\]
together with the H\"older inequality and the Sobolev inequality. However, this results in the exponential-in-time growth of $\|L(t,\cd)\|_{L^2}$, and so is $\|\O(t,\cd)\|_{L^2}$, after applying the Gr\"onwall inequality. Clearly, this is not in line with our goals.
\end{remark}

\qed

Moreover, by H\"older's inequality and Lemma \ref{LEMET11}, one deduces
\[
\begin{split}
|J_2|+|J_3|&\leq\left(\|\tilde{\mathcal{L}}\left(\rho v_r\right)(t,\cd)\|_{L^2}+\left\|\mathcal{L}\left(\rho v_r\right)(t,\cd)\right\|_{L^2}\right)\left\|\partial_z L(t,\cd)\right\|_{L^2}\\
&\leq C\|\rho_0\|_{L^3}\|v(t,\cd)\|_{L^6}\|\p_zL(t,\cd)\|_{L^2}\,.
\end{split}
\]
Using Sobolev imbedding and Young's inequality, one arrives
\[
|J_2|+|J_3|\leq\f{1}{4}\|\na L(t,\cd)\|_{L^2}^2+C\|\rho_0\|_{L^3}^2\|\na v(t,\cd)\|_{L^2}^2\,.
\]
By the same technique, one finds $I_2$ in \eqref{E01} satisfies
\[
\begin{aligned}
|I_2|\leq &\,\|\mathcal{L} \rho(t,\cd)\|_{L^3}\|v(t,\cd)\|_{L^6}\|\nabla L(t,\cd)\|_{L^2}\leq\| \rho_0\|_{L^3}\|v(t,\cd)\|_{L^6}\|\nabla L(t,\cd)\|_{L^2} \\
\leq&\f{1}{4}\|\na L(t,\cd)\|_{L^2}^2+C\|\rho_0\|_{L^3}^2\|\na v(t,\cd)\|_{L^2}^2\,.
\end{aligned}
\]
Combining the above estimates and integrating with the temporal variable over $(0,t)$, one concludes
\be\label{EOO1}
\|L(t,\cd)\|_{L^2}^2+\int_0^t\|\na L(s,\cd)\|_{L^2}^2ds\leq\|L_0\|_{L^2}^2+C\|\rho_0\|_{L^3}^2\int_0^t\|\na v(s,\cd)\|_{L^2}^2ds\leq C_0\,.
\ee
Noting that $L=\O-\mathcal{L}\rho$, this indicates that
\be\label{Fund1}
\begin{split}
\|\O(t,\cd)\|_{L^2}&\leq\|L(t,\cd)\|_{L^2}+\|\mL\rho(t,\cd)\|_{L^2}\leq\|L(t,\cd)\|_{L^2}+\|\rho_0\|_{L^2}\leq C_0\,.
\end{split}
\ee
This finishes the proof.

\qed

\begin{remark}
If Condition \ref{COND1} no longer holds, one can still obtain the following linear temporal growth of $\O$ :
\[
\|\O(t,\cd)\|_{L^2}\leq C_0(1+t)\,.
\]
Readers can find an outline proof in Section \ref{SEC35} of the current paper. This improves related result in \cite[Proposition 4.2]{HR:2010AIHP}, where the authors obtained an exponential-in-time upper bound:
\[
\|\O(t,\cd)\|_{L^2}\leq C_0e^{C_0t}\,.
\]
\end{remark}

\qed

\subsection{$\bl{L^{2^n}}$-bound of $\bl{\o_\th}$}\label{SEC33}
Now we are ready for proving $L^{2^n}$ estimates of $\o_\th$ for $n\in\mathbb{N}$ by the method of induction. Once the following proposition is concluded, Theorem \ref{Main1} is concluded by the Biot-Savart law and interpolations between Lebesgue spaces.
\begin{proposition}\label{PROP3.3}
Let the condition of Proposition \ref{PROP3.2} be satisfied. Suppose $(\o_\th)_0\in L^{2^n}$, with $2\leq n\in\mathbb{N}$. Then
\be\label{EOMEG}
\begin{split}
\left\|\o_\theta(t,\cd)\right\|_{L^{2^n}}^{2^n}+&\int_0^t\big\|\nabla (\o_\theta^{2^{n-1}})(s,\cd)\big\|_{L^2}^2ds+\int_0^t\Big\|\frac{\o_\theta^{2^{n-1}}}{r}(s,\cd)\Big\|_{L^2}^2ds\leq C_{0,n}\,.
\end{split}
\ee
\end{proposition}

\qed

Before the main proof, let us introduce the following lemma on the interpolation:
\begin{lemma}
Given $f\in C_c^\i(\mR^3)$, for any $0<\ve<1$ and $n\in\mathbb{N}$ that $n>2$, the following inequalities hold
\be\label{INT1}
\|f\|^{2^n}_{L^{\f{3\cd 2^n}{2^{n-1}+1}}}\leq \ve\|\na f\|_{L^{3\cd 2^n}}^{2^n}+C_{\ve,n}\|f\|_{L^{2}}^{{c}_{*}}\|\na f\|_{L^2}^2\,;
\ee
\be\label{INT2}
\|f\|^{2^n-2}_{L^{2^n-2}}\leq \ve\|\na (f^{2^{n-1}})\|_{L^{2}}^{2}+C_{\ve,n}\|f\|_{L^{2^{n-1}}}^{\tilde{c}_{*}}\|\na (f^{2^{n-3}})\|_{L^2}^2\,.
\ee
Here ${c}_{*}=2^n-2$, $\tilde{c}_{*}=\f{3\cdot2^{2n-2}}{2^n+6}$, and $C_{\ve,n}$ is a positive constant depending only on $\ve$ and $n$.
\end{lemma}
\pf
Noting that
\[
2<\f{3\cd 2^n}{2^{n-1}+1}<6,\q\forall n\in\mathbb{N},
\]
using Lemma \ref{LEMGN}, one has
\be\label{INT3}
\|f\|_{L^{\f{3\cd 2^n}{2^{n-1}+1}}}\lesssim_n\left\{
\begin{aligned}
&\|f\|_{L^2}^{1-\beta}\|\na f\|_{L^2}^\beta\,;\\
&\|f\|_{L^2}^{1-\gamma}\|\na f\|_{L^{3\cd 2^n}}^\gamma\,.
\end{aligned}
\right.
\ee
Here $\beta=1-2^{-n}$ and $\gamma=\f{2^n-1}{3\cdot2^{n-1}-1}$ are fixed constants. Let $\al\in[0,1]$ be determined later, \eqref{INT3} and Young's inequality indicate
\[
\begin{split}
\|f\|^{2^n}_{L^{\f{3\cd 2^n}{2^{n-1}+1}}}&\leq C_n\left(\|f\|_{L^2}^{1-\beta}\|\na f\|_{L^2}^\beta\right)^{2^n\al}\left(\|f\|_{L^2}^{1-\gamma}\|\na f\|_{L^{3\cd 2^n}}^\gamma\right)^{2^n(1-\al)}\\
&\leq\ve\|\na f\|_{L^{3\cd 2^n}}^{2^n}+C_{\ve,n}\|f\|_{L^{2}}^{c_*}\|\na f\|_{L^2}^{g(\al)}\,,
\end{split}
\]
where
\[
g(\al):=\f{2^n\al\beta}{1-\gamma(1-\al)}.
\]
By the solving the equation $g(\al)=2$, one has $\al=\f{2^{n+1}}{3(2+2^{2n}-3\cdot2^n)}$, which indicates
\[
c_*=2^n-2\,.
\]
 This finishes the proof of \eqref{INT1}. For \eqref{INT2}, by interpolations between Lebesgue spaces, one has
\be\label{INT4}
\|f\|_{L^{2^n-2}}\lesssim_n\left\{
\begin{aligned}
&\|f\|_{L^{3\cd 2^{n-2}}}^{1-\tilde{\beta}}\|f\|_{L^{3\cd 2^{n}}}^{\tilde{\beta}}\,;\\
&\|f\|_{L^{2^{n-1}}}^{1-\tilde{\gamma}}\|f\|_{L^{3\cd 2^n}}^{\tilde{\gamma}}\,.
\end{aligned}
\right.
\ee
Here $\tilde{\beta}=\f{1}{3}-\f{1}{2^{n-1}-1}<\f{1}{3}$, and $\tilde{\gamma}=\f{3(2^n-4)}{5(2^n-2)}$ are fixed constants. Given $\tilde{\al}\in[0,1]$ be determined later, by Young's inequality and \eqref{INT4}, one deduces
\[
\begin{split}
\|f\|^{2^n-2}_{L^{2^n-2}}&\leq\left(\|f\|_{L^{3\cd 2^{n-2}}}^{1-\tilde{\beta}}\|f\|_{L^{3\cd 2^{n}}}^{\tilde{\beta}}\right)^{(2^n-2)\tilde{\al}}\left(\|f\|_{L^{2^{n-1}}}^{1-\tilde{\gamma}}\|f\|_{L^{3\cd 2^n}}^{\tilde{\gamma}}\right)^{(2^n-2)(1-\tilde{\al})}\\
&\leq\ve\|f\|_{L^{3\cd 2^n}}^{2^n}+C_{\ve,n}\|f\|_{L^{2^{n-1}}}^{\tilde{c}_*}\|f\|_{L^{3\cd 2^{n-2}}}^{\tilde{g}(\tilde{\al})}.
\end{split}
\]
Here
\[
\tilde{g}(\tilde{\al}):=\frac{\tilde{\alpha}(1-\tilde{\beta})\left(2^n-2\right) \cdot 2^{n-1}}{2^{n-1}-\left(2^{n-1}-1\right)(\tilde{\alpha}\tilde{\beta}+(1-\tilde{\alpha})\tilde{\gamma})}.
\]
Solving the equation $\tilde{g}(\tilde{\al})=2^{n-2}$, one deduces $\tilde{\al}=\f{2^n+6}{6(2^n+1)}$. This follows that
\[
\tilde{c}_*=\f{3\cdot2^{2n-2}}{2^n+6}\,.
\]
 Thus we conclude \eqref{INT2} by applying the Sobolev inequality.

\qed

\noindent {\bf Proof of Proposition \ref{PROP3.3}. \hspace{1mm}} Now we show the $L^{2^n}$ bound of $\o_\th$. The proof is carried out by an induction argument. First when $n=2$, by taking the $L^4$-energy estimate of the $\o_\th$-equation:
\[
\partial_t \omega_\th+\left(v_r\p_r+v_z\p_z\right)\omega_\th=\left(\Dl-\f{1}{r^2}\right)\o_\th+\frac{v_r}{r} \omega_\th-\p_r\rho\, ,
\]
we arrive at
\be\label{EK0}
\begin{aligned}
\frac{1}{4} \frac{d}{d t}\left\|\o_\theta(t,\cd)\right\|_{L^4}^4+\f{3}{4}\left\|\nabla (\o_\theta^2)(t,\cd)\right\|_{L^2}^2+\Big\|\frac{\o_\theta^2}{r}(t,\cd)\Big\|_{L^2}^2=\un{\int_{\mathbb{R}^3} \frac{v_r}{r}\o_\theta^4 d x}_{K_1}-\un{\int_{\mathbb{R}^3} \partial_r \rho\o_\theta^3 d x}_{K_2}\,.
\end{aligned}
\ee
Using H\"older's inequality, Sobolev imbedding and Young's inequality, one derives
\be\label{EK11}
\begin{split}
|K_1|&\leq\int_{\mR^3}|\O|\,|\o_\th|^3|v_r|dx\leq\|\O(t,\cd)\|_{L^2}\|\o_\th^{2}(t,\cd)\|^{3/2}_{L^6}\|v(t,\cd)\|_{L^4}\\
&\leq\f{1}{8}\|\na(\o_\th^2)(t,\cd)\|_{L^2}^2+C\|\O(t,\cd)\|_{L^2}^4\|v(t,\cd)\|_{L^4}^4\,.
\end{split}
\ee
Using Lemma \ref{LEMGN}, combining interpolations
\[
\|v(t,\cd)\|_{L^4}\lesssim \|v(t,\cd)\|_{L^2}^{2/3}\|\na v(t,\cd)\|_{L^{12}}^{1/3},
\]
and
\[
\|v(t,\cd)\|_{L^4}\lesssim \|v(t,\cd)\|_{L^2}^{1/4}\|\na v(t,\cd)\|_{L^{2}}^{3/4},
\]
one arrives at
\be\label{Ev4}
\begin{split}
\|v(t,\cd)\|_{L^4}^4&\leq C\left(\|v(t,\cd)\|_{L^2}^{8/7}\|\o_\th(t,\cd)\|_{L^{12}}^{4/7}\right)\left(\|v(t,\cd)\|_{L^2}^{4/7}\|\na v(t,\cd)\|_{L^2}^{12/7}\right)\\
&\leq C\|\na (\o_\th^2)(t,\cd)\|_{L^2}^{2/7}\|v(t,\cd)\|_{L^2}^{12/7}\|\na v(t,\cd)\|_{L^2}^{12/7}\,.
\end{split}
\ee
Here we have applied the Biot-Savart law while deriving the first inequality. Substituting \eqref{Ev4} in the far right of \eqref{EK11}, one deduces that
\be\label{EEK1}
\begin{split}
|K_1|&\leq \f{1}{4}\|\na(\o_\th^2)(t,\cd)\|_{L^2}^2+C\|\O(t,\cd)\|_{L^2}^{14/3}\|v(t,\cd)\|_{L^2}^2\|\na v(t,\cd)\|_{L^2}^2\,.
\end{split}
\ee

For term $K_2$, we integrate by parts in the cylindrical coordinates as
\[
\begin{aligned}
K_2&=2\pi \int_{\mathbb{R}} \int_0^{\infty} \partial_r \rho\o_\theta^3 r d r d z \\
&=-2 \pi \int_{\mathbb{R}} \int_0^{\infty} \rho \partial_r\left(\o_\theta^3 r\right) d r d z \\
&=-6 \pi \int_{\mathbb{R}} \int_0^{\infty} \rho\o_\theta^2 \partial_r \o_\theta r d r d z-2\pi\int_{\mathbb{R}^3} \rho \frac{\o_\theta^3}{r} d x\,.
\end{aligned}
\]
Then the Sobolev imbedding theorem, together with H\"older's inequality, Young's inequality and \eqref{ERHO}, indicates that
\be\label{EEK2}
\begin{split}
|K_2|&\leq C\|\rho(t,\cd)\|_{L^{\infty}}\left\|\o_\theta(t,\cd)\right\|_{L^2}\left(\left\|\nabla(\o_\theta^2)(t,\cd)\right\|_{L^2}+\Big\|\frac{\o_\theta^2}{r}(t,\cd)\Big\|_{L^2}\right) \\
&\leq C\|\rho_0\|_{L^{\infty}}^2\left\|\na v(t,\cd)\right\|_{L^2}^2+\frac{1}{4}\left\|\nabla(\o_\theta^2)(t,\cd)\right\|_{L^2}^2+\frac{1}{4}\Big\|\frac{\o_\theta^2}{r}(t,\cd)\Big\|_{L^2}^2 .
\end{split}
\ee
Substituting \eqref{EEK1} and \eqref{EEK2} in \eqref{EK0}, and integrating with the temporal variable over $(0,t)$, one derives
\[
\begin{split}
\left\|\omega_\theta(t, \cdot)\right\|_{L^4}^4&+\int_0^t\left\|\nabla (\omega_\theta^2)(s, \cdot)\right\|_{L^2}^2ds+\int_0^t\Big\|\frac{\omega_\theta^2}{r}(s, \cdot)\Big\|_{L^2}^2ds\\[1mm]
&\leq C\sup_{0\leq s\leq t}\left(\|\O(s,\cd)\|_{L^2}^{14/3}\|v(s,\cd)\|_{L^2}^2\right)\int_0^t\|\na v(s,\cd)\|_{L^2}^2ds+C\|\rho_0\|_{L^{\infty}}^2\int_0^t\left\|\na v(s,\cd)\right\|_{L^2}^2ds\\
&\leq C_0\,.
\end{split}
\]
Here the last inequality follows from \eqref{Fund} and \eqref{Fund1} in previous subsections. This finishes the proof of $n=2$. Now we proceed with the induction by assuming
\be\label{ASM}
\left\|\omega_\theta(t, \cdot)\right\|_{L^{2^k}}^{2^k}+\int_0^t\big\|\nabla (\omega_\theta^{2^{k-1}})(s, \cdot)\big\|_{L^2}^2ds+\int_0^t\Big\|\frac{\omega_\theta^{2^{k-1}}}{r}(s, \cdot)\Big\|_{L^2}^2ds\leq C_{0,k}
\ee
holds for $k=2,3,...,(n-1)$, and we perform the $L^{2^n}$ energy estimate of \eqref{omega} that
\be\label{EK1}
\begin{aligned}
\frac{1}{2^n} \frac{d}{d t}\left\|\o_\theta(t,\cd)\right\|_{L^{2^n}}^{2^n}+\f{2^n-1}{2^{2n-2}}\big\|\nabla (\o_\theta^{2^{n-1}})(t,\cd)\big\|_{L^2}^2+\Big\|\frac{\o_\theta^{2^{n-1}}}{r}(t,\cd)\Big\|_{L^2}^2=\un{\int_{\mathbb{R}^3} \frac{v_r}{r}\o_\theta^{2^n} d x}_{M_1}-\un{\int_{\mathbb{R}^3} \partial_r \rho\o_\theta^{2^n-1} d x}_{M_2}\,.
\end{aligned}
\ee
Here, noticing that
\[
\f{2^{n-1}+1}{3\cd 2^n}+\f{1}{2}+\f{2^n-1}{3\cd 2^n}=1,
\]
one deduces by H\"older's inequality that
\be\label{EK2}
|M_1|\leq\int_{\mR^3}|v_r|\,|\O|\,|\o_\th^{2^n-1}|dx\leq\|v(t,\cd)\|_{L^\f{3\cd 2^n}{2^{n-1}+1}}\|\O(t,\cd)\|_{L^2}\|\o_\th^{2^n-1}(t,\cd)\|_{L^\f{3\cd 2^n}{2^n-1}}.
\ee
Using the Sobolev imbedding, one notices that
\be\label{EK3}
\|\o_\th^{2^n-1}(t,\cd)\|_{L^\f{3\cd 2^n}{2^n-1}}=\|\o_\th^{2^{n-1}}(t,\cd)\|_{L^{6}}^\f{2^n-1}{2^{n-1}}\lesssim_n\|\na( \o_\th^{2^{n-1}})(t,\cd)\|_{L^2}^\f{2^n-1}{2^{n-1}}.
\ee
Substituting \eqref{EK3} into \eqref{EK2} and applying the Young inequality, one deduces
\be\label{EM11}
\begin{split}
|M_1|&\leq\f{2^n-1}{2^{2n+1}}\big\|\nabla (\o_\theta^{2^{n-1}})(t,\cd)\big\|_{L^2}^2+C_n{\|v(t,\cd)\|^{2^n}_{L^\f{3\cd 2^n}{2^{n-1}+1}}\|\O(t,\cd)\|_{L^2}^{2^n}}\\
&\leq\f{2^n-1}{2^{2n+1}}\big\|\nabla (\o_\theta^{2^{n-1}})(t,\cd)\big\|_{L^2}^2+\tilde{C}_{0,n}{\|v(t,\cd)\|^{2^n}_{L^\f{3\cd 2^n}{2^{n-1}+1}}}\,.\\
\end{split}
\ee
Here the last inequality is derived by help of \eqref{Fund1}. Applying \eqref{INT1} and the Biot-Savart law, one derives
\[
\begin{split}
\|v(t,\cd)\|^{2^n}_{L^\f{3\cd 2^n}{2^{n-1}+1}}&\leq\f{2^n-1}{2^{2n+1}C_*\tilde{C}_{0,n}}\|\o_\th(t,\cd)\|_{L^{3\cd 2^n}}^{2^n}+C_{n}\|v(t,\cd)\|_{L^{2}}^{{c}_{*}}\|\na v(t,\cd)\|_{L^2}^2\\
&\leq\f{2^n-1}{2^{2n+1}\tilde{C}_{0,n}}\big\|\nabla (\o_\theta^{2^{n-1}})(t,\cd)\big\|_{L^2}^2+C_{n}\|v(t,\cd)\|_{L^{2}}^{{c}_{*}}\|\na v(t,\cd)\|_{L^2}^2\,.
\end{split}
\]
Here $C_*$ is the best constant of the Sobolev inequality:
\[
\|\na f\|_{L^6}\leq C_*\|\na f\|_{L^2},\q\forall\,f\in H^1.
\]
Inserting this in the far right of \eqref{EM11}, and using the fundamental energy bound \eqref{Fund}, one arrives at
\be\label{EM111}
|M_1|\leq\f{2^n-1}{2^{2n}}\big\|\nabla (\o_\theta^{2^{n-1}})(t,\cd)\big\|_{L^2}^2+C_{0,n}\|\na v(t,\cd)\|_{L^2}^2\,.
\ee
Now it remains to proceed with $M_2$ in \eqref{EK1}. Applying integration by parts, one deduces
\[
M_{2}=-\f{2^n-1}{2^{n-1}}\int_{\mR^3}\rho\o_\th^{2^{n-1}-1}\p_r\o_\th^{2^{n-1}}dx-\int_{\mR^3}\rho\f{\o_\th^{2^{n}-1}}{r}dx\,.
\]
Using Young's inequality and H\"older's inequality, one deduces
\be\label{M21}
|M_2|\leq\f{2^n-1}{2^{2n+1}}\big\|\nabla (\o_\theta^{2^{n-1}})(t,\cd)\big\|_{L^2}^2+\f{1}{2}\Big\|\frac{\o_\theta^{2^{n-1}}}{r}(t,\cd)\Big\|_{L^2}^2+C_n\|\rho_0\|_{L^\i}^2{\int_{\mR^3}|\o_\th|^{2^n-2}dx}\,.
\ee
Utilizing \eqref{INT2}, one has
\[
\int_{\mR^3}|\o_\th|^{2^n-2}dx\leq\f{2^n-1}{2^{2n+1}}\big\|\na (\o_\th^{2^{n-1}})(t,\cd)\big\|_{L^{2}}^{2}+C_{n}\|\o_\th(t,\cd)\|_{L^{2^{n-1}}}^{\tilde{c}_{*}}\big\|\na (\o_\th^{2^{n-3}})(t,\cd)\big\|_{L^2}^2\,.
\]
Inserting this in the far right of \eqref{M21}, and using the assumption \eqref{ASM}, one arrives at
\be\label{EM222}
|M_2|\leq\f{2^n-1}{2^{2n}}\big\|\nabla (\o_\theta^{2^{n-1}})(t,\cd)\big\|_{L^2}^2+\f{1}{2}\Big\|\frac{\o_\theta^{2^{n-1}}}{r}(t,\cd)\Big\|_{L^2}^2+C_{0,n}\big\|\na (\o_\th^{2^{n-3}})(t,\cd)\big\|_{L^2}^2\,.
\ee
Substituting \eqref{EM111} and \eqref{EM222} in \eqref{EK1}, then integrating with the temporal variable over $(0,t)$, one deduces
\[
\begin{split}
\left\|\o_\theta(t,\cd)\right\|_{L^{2^n}}^{2^n}+&\int_0^t\big\|\nabla (\o_\theta^{2^{n-1}})(s,\cd)\big\|_{L^2}^2ds+\int_0^t\Big\|\frac{\o_\theta^{2^{n-1}}}{r}(s,\cd)\Big\|_{L^2}^2ds\\
&\leq C_{0,n}\int_0^t\left(\|\na v(s,\cd)\|_{L^2}^2+\big\|\na (\o_\th^{2^{n-3}})(t,\cd)\big\|_{L^2}^2\right)ds\leq C_{0,n}\,.
\end{split}
\]
Here the last inequality holds due to the fundamental energy bound \eqref{Fund} and the assumption of induction \eqref{ASM}. This finishes the proof.

\qed

\subsection{Higher-order estimates}\label{SEC34}
We carry out the proof of Theorem \ref{Main2} in this subsection. To do this, we first \textbf{claim} that
\be\label{ENU}
\int_0^{t}\|\nabla v(s,\cdot)\|_{L^\infty}ds\lesssim t^{\left(4/5\right)_+}\,.
\ee
Suppose \eqref{ENU} is achieved. By acting $\nabla$ on $\eqref{Bous}_2$, we know that
\[
\p_t\nabla\rho+v\cdot\nabla\nabla\rho=-\nabla v\cdot\nabla\rho.
\]
The routine $L^p$ ($1\leq p\leq\i$) estimate follows that
\[
\|\nabla\rho(t,\cdot)\|_{L^p}\leq\|\nabla\rho_0\|_{L^p}+\int_0^t\|\nabla v(s,\cdot)\|_{L^\infty}\|\nabla\rho(s,\cdot)\|_{L^p}ds.
\]
By the Gr\"onwall inequality and using \eqref{ENU}, we arrive
\[
\sup_{0\leq s\leq t}\|\nabla\rho(s,\cdot)\|_{L^p}\leq\|\nabla\rho_0\|_{L^p}\exp\left(\int_0^{t}\|\nabla v(s,\cdot)\|_{L^\infty}ds\right)\leq C_{0}\exp\left(t^{\left(4/5\right)_+}\right)\,.
\]
This concludes \eqref{EH1} in Theorem \ref{Main2}. Apply $\na^m$ $(m\in\mathbb{N},\,\,m\geq3)$ to \eqref{Bous}$_{1,2}$ to derive that
\begin{equation}\label{MB3}
\left\{
\begin{aligned}
&
\p_t\na^mv+v\cdot\nabla\na^mv+\nabla\na^mP-\Delta\na^mv=\na^m(\rho \bl{e_3})-[\na^m,v\cdot\nabla]v,\\[4mm]
&\p_t\na^m\rho+v\cdot\na\na^m\rho=-[\na^m,v\cdot\nabla]\rho.\\
\end{aligned}
\right.
\end{equation}
Performing the $L^2$ energy estimate of \eqref{MB3}, we have
\[
\bali
\frac{1}{2}\frac{d}{dt}\left\|\na^m (v,\rho)(t,\cdot)\right\|_{L^2}^2+\left\|\na^{m+1}v(t,\cdot)\right\|_{L^2}^2=&-\int_{\mathbb{R}^3}[\na^m,v\cdot\nabla]v\na^mvdx-\int_{\mathbb{R}^3}[\na^m,v\cdot\nabla]\rho\na^m\rho dx\\
&+\int_{\mathbb{R}^3}\na^m(\rho \bl{e_3})\na^mvdx.
\eali
\]
By Lemma \ref{LEMET1} and the Cauchy-Schwarz inequality, the above equation implies
\[
\bali
&\frac{d}{dt}\left\|\na^m (v,\rho)(t,\cdot)\right\|_{L^2}^2+\left\|\na^{m+1}v(t,\cdot)\right\|_{L^2}^2\lesssim\|\nabla^m(v,\rho)(t,\cdot)\|_{L^2}^2\left(\|\nabla(v,\rho)(t,\cdot)\|_{L^\infty}+1\right).
\eali
\]
Using the Gr\"onwall inequality, one deduces that
\[
\begin{split}
\left\|\na^m (v,\rho)(t,\cdot)\right\|_{L^2}^2\leq&C_0\|\na^m(v_0,\,\rho_0)\|_{L^2}^2\exp\left( C_{0}\int^t_{0}\left(1+\|(\na v,\nabla\rho)(s,\cdot)\|_{L^\infty}\right)ds\right)\\
\les_{0,m}&\exp\left(\exp\left(t^{\left(4/5\right)_+}\right)\right),\quad \forall t\in[0,\i)\,.
\end{split}
\]
This concludes \eqref{EH2} in Theorem \ref{Main2}.

\qed

\noindent {\bf Proof of the equation \eqref{ENU}. \hspace{1mm}} Denoting $\o=\na\times v$ and acting $\na\,\times$ on \eqref{Bous}$_{1}$, one deduces that
\bes
\lt\{
\begin{aligned}
&\p_t \o-\Delta \o=\nabla\times(v\cdot\nabla v)+\nabla\times(\rho \bl{e_3});\\[4mm]
&\o(0,x)=\nabla\times v_0(x).
\end{aligned}
\rt.
\ees
For the further convenience, we split $\o$ into two parts:
\[
\o:=\o_1+\o_2,
\]
where $\o_1$ solves the linear parabolic equation with the initial value $\nabla\times v_0(x)$:
\bes
\lt\{
\begin{aligned}
&\p_t \o_1-\Delta \o_1=0;\\[4mm]
&\o(0,x)=\nabla\times v_0(x).
\end{aligned}
\rt.\ees
Clearly, $\o_1$ is smooth when $t$ is strictly away from zero, and it does not grow as $t\to\i$. So we only consider the rest part. Noticing that $\o_2$, which has homogeneous initial data, satisfies
\bes
\p_t \o_2-\Delta \o_2=\nabla\times(v\cdot\nabla v)+\nabla\times(\rho \bl{e_3}).
\ees
Using the estimate \eqref{EOMEG}, Biot-Savart law, and the Sobolev imbedding theorem, one deduces
\[
\sup_{0\leq s\leq t}\|(v\cdot\na v)(s,\cd)\|_{L^{2^n}}\leq\sup_{0\leq s\leq t}\left(\|v(s,\cd)\|_{L^{\i}}\|\na v(s,\cd)\|_{L^{2^n}}\right)\leq C_{0,n}\,.
\]
This implies
\[
\left(\int_0^t\|(v\cdot\na v)(s,\cd)\|^p_{L^{2^n}}ds\right)^{1/p}\leq C_{0,n}t^{1/p},\q\forall p\in (1,\i)\,.
\]
Meanwhile, by \eqref{ERHO}, it is clear that
\[
\left(\int_0^t\|\rho(s,\cd)\|^p_{L^{2^n}}ds\right)^{1/p}\leq\|\rho_0\|_{L^{2^n}}\left(\int_0^tds\right)^{1/p}\leq C_{0}t^{1/p},\q\forall p\in(1,\i)\,.
\]
Thus, by applying the maximal regularity of the heat flow (Lemma \ref{MRP}) and the Biot-Savart law, it is clear that
\be\label{Ew2}
\left(\int_0^t\|\na^2v(s,\cd)\|^p_{L^{2^n}}ds\right)^{1/p}\leq C_{0,p,n}t^{1/p}\,.
\ee
Then utilizing Lemma \ref{LEMGN}  and the H\"older inequality, we find
\[
\begin{split}
\int_0^{t}\|\nabla v(s,\cdot)\|_{L^\infty}ds&\leq C_{0,n}\,\int_0^t\|\nabla v(s,\cd)\|_{L^2}^{1-\zeta_n}\|\nabla^2 v(s, \cdot)\|_{L^{2^n}}^{\zeta_n} d s\\
&\leq C_{0,n}\,\left(\int_0^t\|\na v(s,\cd)\|^2_{L^2}ds\right)^{\f{1-\zeta_n}{2}}\left(\int_0^t\|\na^2v(s,\cd)\|^{\f{2\zeta_n}{1+\zeta_n}}_{L^{2^n}}ds\right)^{\f{1+\zeta_n}{2}}\\
&\leq C_{0,n}\,t^{\f{1+\zeta_n}{2}}\,,
\end{split}
\]
where
\[
\zeta_n=\f{1}{2}\left(\f{5}{6}-\f{1}{2^n}\right)^{-1}.
\]
Here the last inequality follows from \eqref{Ew2} and \eqref{Ev}. Noting that
\[
\zeta_n\searrow \f{3}{5}\q\text{as}\q n\to\i\,,
\]
thus we conclude \eqref{ENU} by choosing $n$ sufficiently large. This finishes the proof of the claim.

\subsection{On Proposition \ref{COR1} and Proposition \ref{Main3}}\label{SEC35}
When Condition \ref{COND1} is discarded, results in Theorem \ref{Main1} and \ref{Main2} turn to Proposition \ref{COR1} and Proposition \ref{Main3}, respectively. With the help of the proof carried out in previous sections, we outline the proofs for these two propositions below.

Proceeding the standard energy estimate of \eqref{Bous}, we have
\be\label{C00}
\f{1}{2}\f{d}{dt}\|v(t,\cd)\|_{L^2}^2+\|\na v(t,\cd)\|_{L^2}^2\leq \|v(t,\cd)\|_{L^2}\|\rho(t,\cd)\|_{L^2}\,.
\ee
Cancelling $\|v(t,\cd)\|_{L^2}$ on each side, one arrives
\[
\f{d}{dt}\|v(t,\cd)\|_{L^2}\leq \|\rho(t,\cd)\|_{L^2}\,.
\]
Integrating on both sides and applying \eqref{ERHO}, one deduces
\be\label{C01}
\|v(t,\cd)\|_{L^2}\leq \|v_0\|_{L^2}+t\|\rho_0\|_{L^2}.
\ee
Substituting \eqref{C01} in \eqref{C00} and then integrating with the temporal variable, one derives
\be\label{C02}
\|v(t,\cd)\|_{L^2}^2+\int_0^t\|\na v(s,\cd)\|_{L^2}^2ds\leq C_0(1+t)^2\,.
\ee
Therefore, estimate \eqref{EOO1} degenerates to
\[
\|L(t,\cd)\|_{L^2}^2+2\int_0^t\|\na L(s,\cd)\|_{L^2}^2ds\leq \|L_0\|_{L^2}^2+ C\|\rho_0\|_{L^3}^2\int_0^t\|\na v(s,\cd)\|_{L^2}^2ds\leq C_0(1+t)^2\,,
\]
and thus
\[
\|\O(t,\cd)\|_{L^2}\leq C_0(1+t)\,.
\]
Following the procedure in Section \ref{SEC33}, one has
\[
\left\|\o_\theta(t,\cd)\right\|_{L^{2^n}}^{2^n}+\int_0^t\big\|\nabla (\o_\theta^{2^{n-1}})(s,\cd)\big\|_{L^2}^2ds+\int_0^t\Big\|\frac{\o_\theta^{2^{n-1}}}{r}(s,\cd)\Big\|_{L^2}^2ds\leq C_{0,n}(1+t)^{M_n}\,,\q\text{for }n=2,3,4,...\,.
\]
In this way, using the maximal regularity of the heat flow, one rewrites \eqref{Ew2} as
\be\label{EMM}
\left(\int_0^t\|\na^2v(s,\cd)\|^p_{L^{2^n}}ds\right)^{1/p}\leq C_{0,p,n}t^{M}\,.
\ee
This proves Proposition \ref{COR1}. Interpolating \eqref{C02} and \eqref{EMM}, one has
\[
\int_0^{t}\|\nabla v(s,\cdot)\|_{L^\infty}ds
\]
grows at most algebraically as $t\to\i$. Following the higher-order energy estimate in Section \ref{SEC34}, one concludes Proposition \ref{Main3}.

\qed

\section*{Acknowledgments}
\addcontentsline{toc}{section}{Acknowledgments}
\q\ The author wishes to thank the anonymous referee for various comments which have improved this manuscript. He also would like to thank Prof. Xin Yang, Dr. Chulan Zeng in UC Riverside, and Prof. Xinghong Pan in Nanjing University of Aeronautics and Astronautics for their helpful discussions on the current paper.

Z. Li is supported by National Natural Science Foundation of China (No. 12001285) and Natural Science Foundation of Jiangsu Province (No. BK20200803).

\medskip
\medskip

{\footnotesize

{\sc Z. Li: School of Mathematics and Statistics, Nanjing University of Information Science and Technology, Nanjing 210044, China}

  {\it E-mail address:}  zijinli@nuist.edu.cn

}
\end{document}